\documentclass[a4paper]{amsart}
\pdfoutput=1
\usepackage[utf8]{inputenc}
\usepackage[english]{babel}
\usepackage[T1]{fontenc}
\usepackage{amsmath}
\usepackage{hyperref}

\usepackage{amssymb,color}
\usepackage{amsfonts}
\usepackage{amsmath}
\usepackage{euscript}
\usepackage{enumerate}
\usepackage{graphics}
\usepackage{tikz}
\usetikzlibrary{arrows}

\newtheorem{theorem}{Theorem}[section]
\newtheorem{proposition}[theorem]{Proposition}

\theoremstyle{definition}
\newtheorem{definition}[theorem]{Definition}

\numberwithin{equation}{section}

\begin{document}

\title{Branching laws for square integrable representations}
\date{today}
\author{Jorge A. Vargas}
\email{vargas@famaf.unc.edu.ar}
\thanks{CONICET, FAMAF, Aarhus University }.

\begin{abstract} We present an overview of results on branching laws for square integrable representations of a semisimple Lie group,  restricted to a closed reductive subgroup. The overview is partial and it is based on joint work with Bent {\O}rsted and the deep work of Toshiyuki Kobayashi.
\end{abstract}

\maketitle

 \section{Introduction} Let $H$ denote a connected reductive Lie group and  $(\pi,V)$ a unitary representation for $H$. To follow we sketch the decomposition of $\pi$ as a direct integral. The result is: \begin{equation}\label{eq:direct} V\equiv \int_{i\in \hat H} V_i \otimes M_i\,\, d\mu (i). \end{equation}
 Here, $\hat H$ is the set of equivalence classes of unitary representations for $H$. For each $i \in \hat H $, $V_i$ is the representation space of an irreducible, unitary representation $(\pi_i, V_i)$ for $H$; $M_i$ is a Hilbert space, $\mu$ is a measure on $\hat H$. The right hand side, roughly speaking, is the space of square integrable functions  from $\hat H$ into the disjoint union $\bigsqcup_{i \in \hat H},  V_i \otimes M_i$ so that for each $i$  the value of the function on $i\in \hat H$ belongs to $V_i \otimes M_i$.
 The equivalence is described as follows: We recall a vector $v \in V$ is {\it smooth} if the function $H\ni h \rightarrow \pi(h)v$ is smooth. The subspace of smooth vectors $V^\infty$ is endowed with its usual topology. Then, there exists a family of $H$ equivariant linear maps $P_i : V^\infty \rightarrow V_i^\infty \otimes M_i$, continuous in smooth topology, so that \begin{equation}\label{eq:norm} \Vert v \Vert^2 =\int_{i \in \hat H} \Vert P_i(v) \Vert^2 d\mu(i).  \end{equation}
 Two simple examples of the decomposition are:

 {\it Fourier series.} Let $S\equiv \mathbb R/\mathbb Z$ denote the module one complex numbers. We endow $S$ with Lebesgue measure.  Let $V=L^2(S)$. Then, $L_s(f)(z):=f(s^{-1}z)$ defines a unitary representation of $S$ in the Hilbert space $L^2(S)$. For each $n\in \mathbb Z$, $S$ has an irreducible representation in $\mathbb C$ defined by the formula $\pi_n(z)=z^n$. The representation $(\pi_n, \mathbb C)$ is equivalent to the representation $L_{\bullet}$ restricted to the subspace $\mathbb C z^{-n}$. From Fourier analysis, we know \begin{equation} \label{eq:fourier} L^2(S) \ni f =\sum_{n\in \mathbb Z} \hat f(n) z^n. \end{equation} Next, we rewrite the formula \eqref{eq:fourier} in the language of direct integral. Here $\hat H\equiv \mathbb Z$ via the map $(\pi_n, \mathbb C)\leftarrow n $, $M_n=\mathbb C$, the measure $\mu$ is $\mu(\pi_n, \mathbb C)=1$, the function $P_n$ is  $P_n(f)=\hat f(n)$. Therefore \eqref{eq:fourier} becomes $L^2(S)\equiv Hilbert(\sum_{n \in \mathbb Z} (\pi_n, \mathbb C))$.  In this case, subspace of smooth vectors is equal to subspace of smooth functions, and,  since square integrable functions on $S$ are integrable functions, $P_n$ extends to a linear map from $L^2(S)$ into $\mathbb C$.

 \medskip

 {\it Fourier integrals.} In this case $H=\mathbb R$, $V=L^2(\mathbb R)$ the space of square integrable functions in $\mathbb R$ with respect to Lebesgue measure. The representation of $\mathbb R$ in $L^2(\mathbb R)$ is by left translation $L_r(f)(x)=f(x-r), r,x \in \mathbb R, f\in L^2(\mathbb R)$. For each $\lambda \in \mathbb R$ we have an unitary irreducible representation of the group $\mathbb R$ in $\mathbb C$ defined via the formula $\pi_\lambda(x)=exp(i\lambda x)$. Then, the spectral theorem yields there is a equivalence $\hat{ \mathbb R} \ni class((\pi_\lambda, \mathbb C)) \leftarrow \lambda \in \mathbb R$. We set as measure $\mu$ on $\hat{ \mathbb R}$ the image of Lebesgue measure divided by $2\pi$ via the equivalence. The subspace $M_\lambda=\mathbb C$.  The space of smooth vectors in the representation  $V=L^2(\mathbb R)$ is   the subspace of tempered functions, the smooth topology agrees with the Schwartz topology in the space of tempered functions. For each smooth vector (tempered function)$f$ we have the Fourier transform $\hat f(\lambda)=\int_\mathbb R f(x) exp(-i\lambda x) dx =:P_\lambda (f)$. Then, $P_\lambda : V^\infty \rightarrow \mathbb C$ is a continuous $\mathbb R$-equivariant linear map, and from the properties of Fourier integrals we have the equality $\Vert f \Vert_{L^2(\mathbb R)} =\frac{1}{2\pi} \int_{\lambda \in \hat{\mathbb R} } P_\lambda(f)(\lambda) d\lambda$. Therefore, $L^2(\mathbb R)\equiv \int_{\lambda \in \hat{\mathbb R}} (\pi_\lambda, \mathbb C) d\lambda$. It is a simple fact to show that $P_\lambda$ does not extend  to a continuous linear functional for $L^2(\mathbb R)$.

 \subsection{Problems} A problem possed by either Physics or Number Theory or other areas of mathematics is to explicit  \eqref{eq:direct}, \eqref{eq:norm} as we have presented in the two examples developed by Fourier at the beginning of the nineteen century. For the case of an abelian finite group, the problem was solved by Frobenius around 1890. For compact groups the problem was solved by H. Weyl and F. Peter around 1930. The first advance for noncompact groups was in the 1940's due to the work of Bargman, Gelfand I., Harish-Chandra on the left regular representation for  Lorentz groups $SO(2,1), SO(3,1)$. Then, mainly Harish-Chandra, computed \eqref{eq:direct}, \eqref{eq:norm} for the left regular representation of a semisimple Lie group $H$. The work of Harish-Chandra on \eqref{eq:direct}, \eqref{eq:norm} extended for over   thirty years and have important consequences, for the time being, in partial differential equations, spectral theory, number theory \cite{M} \cite{Ve}.  A problem that remains is to carry out the proposed work for the unitarily induced representation $Ind_K^H(\sigma)=L^2(H \times_\sigma W)$ where $H$ is a semisimple Lie group, and $(\sigma, W)$ is a unitary irreducible representation of a noncompact closed  subgroup $K$ of $G$. For example, $L^2(SL(n,\mathbb R)/SL(n,\mathbb Z))$ is far from been understood \cite{M}. Other problem that remains unsolved, on which we will comment in more detail, is to analyze \eqref{eq:direct}, \eqref{eq:norm} for the restriction to $H$ of an irreducible unitary representation of a  group $G$ that    contains $H$ as a closed subgroup.  The work of N.  Wallach, the Dutch school, Japanese School, Danish School,  and the French school is important and all of them have obtained relevant contributions on the problem.

 \subsection{Discrete series}For some semisimple Lie groups $G$, there exist  unitary irreducible representation  with the property that any of its matrix coefficients are square integrable with respect to   Haar measure of $G$. The set of equivalence classes of irreducible square integrable representations of $G$ is appel the Discrete Series for $G$, and any square integrable representation is called a Discrete Series.  Thanks to the work of Harish-Chandra, Schmid, Enright-Wallach, Duflo, Wallach, Connes-Moscovici, Hotta, R. Parthasarathy, now a days, we have a good understanding of such representations. One realization of the space $V$ for a Discrete Series  is as the space of the $L^2$-solutions of an elliptic operator. More precisely, we fix a maximal compact subgroup $K$ for $G$. For example, for the group $ SU(1,1)     $ the subgroup of diagonal matrices is a maximal compact subgrop. A theorem of Harish-Chandra says: an arbitrary semisimple Lie group $G$ admits a  Discrete Series representation if and only if the rank of $K$ is equal to the rank of $G$. Equivalently, if each torus in $K$ is a Cartan subgroup of $G$.  For example,  $SO(2n,1)$ has discrete series representations, whereas, $SO(2n+1,1)$  does not have discrete series representations. We now fix an irreducible representation $(\tau, W)$ of $K$ and consider the fiber bundle $W\rightarrow G \times_{\tau} W \rightarrow G/K$. It has been shown that this bundle carries a $G$-invariant metric and that the homogeneous space $G/K$ carries a $G$-invariant Radon measure. Thus, the space of $L^2$-sections, $L^2(G\times_\tau W)$, of the bundle is well defined. The space $L^2(G\times_\tau W)$ is represented by the space \begin{equation*} L^2(G\times_\tau W)   := \{ f \in L^2(G) \otimes W :   f(gk)=\tau(k)^{-1} f(g),   g  \in G, k \in K \}.
   \end{equation*}

   Here, the action of $G$ is by left translation $L_x, x \in G.$   The inner product on $L^2(G)\otimes W$ is given by \begin{equation*}(f,g)_{L^2(G)} =\int_G (f(x),g(x))_W dx, \end{equation*} where $(...,...)_W$ is a $K-$invariant inner product on $W.$

    Subsequently, $L_D $ (resp. $R_D)$ denotes the left infinitesimal (resp. right infinitesimal) action on functions from $G$ of an element  $D$ in universal enveloping algebra $U(\mathfrak g)$ for the Lie algebra $\mathfrak g$.  As usual, $\Omega_G$ denotes the Casimir operator for $\mathfrak g.$   Following Hotta-Parthasarathy, Enright-Wallach, Atiyah-Schmid, we assume  $rank K=rank G$ and we fix a maximal torus $T$ of $K$.  We set $\rho$ to be one half of the sum of the elements in a system of positive roots for the pair $(\mathfrak g, \mathfrak t)$ and let $\lambda \in Lie(T)$ so that $\lambda +\rho$ is the differential of a character of $T$. Finally, we consider the space
\begin{eqnarray*}
 \lefteqn{V_{\tau,\lambda}^G:= H^2(G, \tau_)\lambda  :=\{ f \in L^2(G) \otimes W : f(gk)=\tau(k)^{-1} f(g)} \hspace{3.0cm} \\ & & g\in G, k \in K, R_{\Omega_G} f= [(\lambda, \lambda) -(\rho, \rho)] f  \}.
 \end{eqnarray*}
We also  recall, $R_{\Omega_G}=L_{\Omega_G} $ is an elliptic $G-$invariant operator on the vector bundle $W \rightarrow G \times_\tau W \rightarrow G/K$ and hence,  $\,H^2(G,\tau)_\lambda$ consists of smooth sections, moreover point evaluation $e_x$ defined by  $ \,H^2(G,\tau)_\lambda \ni f \mapsto f(x) \in W $ is continuous for each $x \in G$. Therefore, the orthogonal projector $P_\tau$ onto $\,H^2(G,\tau)_\lambda$ is an integral map (integral operator) represented by the smooth {\it matrix  kernel} or {\it reproducing kernel}     \begin{equation} \label{eq:Klambda}K_\tau : G\times G \rightarrow End_\mathbb C (W) \end{equation} which satisfies $  K_\tau (\cdot ,x)^\star w$ belongs to $\,H^2(G,\tau)_\lambda$ for each $x \in G, w \in W$ and $$ (P_\tau (f)(x), w)_W=\int_G (f(y), K_\tau (y,x)^\star w)_W dy, \,     f\in L_2(G\times_\tau W).$$ It has been shown by Schmid, Enright-Wallach, Hotta et all.
\begin{theorem}\label{thm:ewhp} Whenever the highest weight of $\tau$ is equal to $\lambda +\rho -2\rho_c$\footnote{ $\rho_c$ is defined as $\rho$ for the pair $(\mathfrak k,\mathfrak t)$} ,   $H^2(G,\tau)_\lambda$ is a irreducible square integrable representation for $G$ and varying $\lambda$ so that $(\lambda, \alpha)\not=0 \,\,\forall \alpha$ in the root system for $(\mathfrak g, \mathfrak t)$,  we obtain the totality of the Discrete Series for $G$. That is, any irreducible square integrable representation for $G$ is equivalent to a unique $H^2(G,\tau)_\lambda$.
\end{theorem}
We want to understand the formula \eqref{eq:direct} for the representation $V_\tau^G$ and as group $H$ we consider a closed reductive subgroup of $G$. For this, we consider the subset  of points in  $i\in \hat H$  so that $\mu(\{i\})>0$. This subset is called the subset of $H$-discrete factors. The subset of $H$-discrete points is equal to the subset of $i\in \hat H$ so that $V_i$ is a subrepresentation of $V$.
It has  been shown
\begin{theorem} \label{thm:disfac} (Gross-Wallach, Vargas) Let $V=V_\tau^G$ be a discrete series for $G$. Fix a closed reductive subgroup $H$ of $G$. Then, for the decomposition \eqref{eq:direct} for $res_H(V)$, all the  $H$-irreducible representations in the support of $\mu$ are tempered representations for $H$. Moreover,  the $H$-discrete factors are again square integrable representations for $H$.
\end{theorem}

 \section{Structure of intertwining maps} We fix $G,K,H, V_\tau^G$ as in the previous section. Thus, Theorem~\ref{thm:disfac} implies whenever a irreducible representation of$H$ is a subrepresentation of the restriction to $H$, $res_H(V_\tau^G)$, of $V_\tau^G$, we have the subrepresentation is an irreducible square integrable representation for $H$.

 To continue, we fix a maximal compact subgroup $L$ of $H$. Thus, Theorem~\ref{thm:ewhp} and Theorem~\ref{thm:disfac} yields each $H$-factor of $res_H(V_{\tau,\lambda}^G )=res_H(H^2(G,\tau)_\lambda)$ is   equivalent to a representation $(L_\bullet^H, H^2(H,\sigma)_\mu)$ for   some $\mu \in Lie(T\cap H)$ as in  Theorem~\ref{thm:ewhp}.
 Theorem~\ref{thm:ewhp} implies that both spaces $H^2(G,\tau)_\lambda$, $H^2(H,\sigma)_\mu$ are solutions spaces of elliptic equations. Therefore point evaluation on each space is continuous in $L^2$-topology. A consequence of this is: \\ {\it Each continuous linear map  either from $H^2(G,\tau)_\lambda$ to  $H^2(H,\sigma)_\mu$ or from , $H^2(H,\sigma)_\mu$ to $H^2(G,\tau)_\lambda$  is an integral operator represented by a Carleman kernel}.

  More precisely, we fix
    a continuous intertwining linear $H-$map  $T : H^2(H\times_\sigma Z)\rightarrow \,H^2(G,\tau)_\lambda$.

  Then, for each $x \in G, w \in W$ the linear functional  $H^2(H\times_\sigma Z) \ni g \mapsto (Tg(x),w)_W$ is continuous.  Riesz representation Theorem shows there exists a function  $$K_T : H\times G \rightarrow Hom_\mathbb C (Z,W)$$ so that
the map $ h\mapsto K_T(h,x)^\star(w)$ belongs  to $H^2(H \times_\sigma Z)$ and for $g \in H^2(H\times_\sigma Z), x \in G,  w \in W$  we have the absolutely  convergent integral and the equality
\begin{equation}\label{eq:tisintegral}
(Tg(x),w)_W    =\int_H (g(h), K_T(h,x)^\star w)_Z  dh.
\end{equation}
That is, $T$  is the integral map  $$Tg(x)=\int_H K_T(h,x) g(h) dh,  x \in G.$$

 $K_T$ is a smooth function,     $ K_T(h,\cdot)z=K_{T^\star}(\cdot,h)^\star z \in H^2(G,\tau)_\lambda$ and $ K_T(e,\cdot)z \in H^2(G,\tau)_\lambda[V_\mu^H][V_{\mu+\rho_n^H}^L]$ is a $L$-finite vector in $  H^2(G,\tau)_\lambda$.

 Similarly, we obtain a proof  for intertwining maps from  $H^2(G,\tau)_\lambda$ into $H^2(H,\sigma)_\mu$.

 \begin{definition} A linear map $S$ from $H^2(G,\tau)_\lambda$ into $H^2(H,\sigma)_\mu$ is represented via a differential operator, if there exists a linear differential operator $D $ from the space of smooth sections of $G\times_\tau W\rightarrow G/K$ into the space of smooth sections of $G\times_\sigma Z \rightarrow G/L$ so that $S(f)(h)= D(f)(h)$, for every $h \in H$ and smooth section $f$ for $G\times_\tau W\rightarrow G/K$.
 \end{definition}

 \begin{proposition}({\O}rsted-Vargas) a)Whenever a linear map $S: H^2(G,\tau)_\lambda \rightarrow H^2(H,\sigma)_\mu$ is represented via a differential operator, then $S$ is continuous in $L^2$-topologies.

 b) A continuous linear map $S: H^2(G,\tau)_\lambda \rightarrow H^2(H,\sigma)_\mu$ is represented via a differential operator if and only if $K_S(\cdot,e)^\star z \in H^2(G,\tau)_\lambda$ is a $K$-finite vector for each $z\in Z$.
 \end{proposition}

 \section{Discretely decomposable representations} A particular decomposition in  \eqref{eq:direct} is when $\mu$ is a discrete measure. That is the support of $\mu$ is at most a countable set and $\mu(\{i\})>0$ for each $i$ in the support. This representations have a special name. Namely,
 \begin{definition} \label{dfn:disdecom} A representation $(\pi , V $ is discretely decomposable over $H$, if there exists an orthogonal family of closed, $H-$invariant, $H-$irreducible subspaces of $V $ so that the closure of its algebraic sum is equal to $V.$   \end{definition}
 When $\dim M_i$ is finite for every $i$, we   have a name for such a representations.
\begin{definition}  A representation $(\pi , V )$ is $H-$admissible if the representation is discretely decomposable and the  multiplicity of each irreducible factor is finite. \end{definition}
In Kobayashi-Oshima, Vargas,   has found a complete list of triples $(G,H, (\pi,V) )$ such that $(G,H)$ is a symmetric pair and $(\pi,V)$ is a discrete series for $G$ with  an $H-$admissible representation.

We have obtained a criterion for either discretely decomposable or admissibility.
\begin{theorem}({\O}rsted-Vargas) Let $(\pi,V=H^2(G,\tau)_\lambda )$ be a discrete series representation for $G$. Then,

a) $res_H(\pi)$ is $H$-discretely decomposable if and only if  there exists a discrete series $H^2(H,\sigma)_\mu$ for $H$ and a nonzero intertwining linear map $S$ from $V$ into $H^2(H,\sigma)_\mu$ represented via a differential operator.

b) If every intertwining linear map from $H^2(G,\tau)_\lambda $  into $H^2(H,\sigma)_\mu$ is represented via a differential operator, then the multiplicity of $H^2(H,\sigma)_\mu$ in $H^2(G,\tau)_\lambda $ is finite.
\end{theorem}
Other criteria for discrete decomposability is obtained by means of the spherical function attached to the lowest $K$-type of $(\pi, H^2(G,\tau)_\lambda)$. Let $\mathfrak z_\mathfrak h$ denote the center of the universal enveloping algebra of $\mathfrak h$. A smooth function $f$on $H$ is called $\mathfrak z_\mathfrak h$-finite if the dimension of the span of $\{L_D (f): D \in \mathfrak z_\mathfrak h  \}$ is finite. Let $P$ denote the orthogonal projector onto the lowest $K$-type $W$ of $H^2(G,\tau)_\lambda$. Let $\Phi (g)=P\pi(g)P$ denote the spherical function attached to the lowest $K$-type $(\tau, W)$ of $\pi$.
\begin{theorem}({\O}rsted-Vargas) $res_H(H^2(G,\tau)_\lambda)$ is $H$-discretely decomposable if and only if $\Phi$ is a $\mathfrak z_\mathfrak h$-finite function.

\end{theorem}

  \bibliographystyle{plain}

\end{document}